\documentclass{amsart}
\usepackage{pinlabel}
\usepackage{amsmath, amsthm, amssymb, amscd, mathrsfs, eucal, epsfig, color}
\usepackage{hyperref}

\input xy
\xyoption{all}

\begin{document}

\newtheorem*{theorem_A}{Theorem A}
\newtheorem{theorem}{Theorem}[section]
\newtheorem{lemma}[theorem]{Lemma}
\newtheorem{corollary}[theorem]{Corollary}
\newtheorem{conjecture}[theorem]{Conjecture}
\newtheorem{proposition}[theorem]{Proposition}
\newtheorem{question}[theorem]{Question}
\newtheorem*{answer}{Answer}
\newtheorem{problem}[theorem]{Problem}
\newtheorem*{main_theorem}{Main Theorem}
\newtheorem*{claim}{Claim}
\newtheorem*{criterion}{Criterion}
\theoremstyle{definition}
\newtheorem{definition}[theorem]{Definition}
\newtheorem{construction}[theorem]{Construction}
\newtheorem{notation}[theorem]{Notation}
\newtheorem{convention}[theorem]{Convention}
\newtheorem*{warning}{Warning}
\newtheorem*{assumption}{Simplifying Assumptions}

\theoremstyle{remark}
\newtheorem{remark}[theorem]{Remark}
\newtheorem{example}[theorem]{Example}
\newtheorem{scholium}[theorem]{Scholium}
\newtheorem*{case}{Case}

\def\id{\text{id}}
\def\H{\mathbb H}
\def\Z{\mathbb Z}
\def\N{\mathbb N}
\def\R{\mathbb R}
\def\RR{\mathcal R}
\def\C{\mathbb C}
\def\Chat{{\hat{\C}}}
\def\D{\mathbb D}
\def\E{\mathbb E}
\def\RR{\mathcal{R}}
\def\SS{\mathcal{S}}
\def\DL{\mathcal{DL}}
\def\P{\mathcal{P}}
\def\M{\mathcal{M}}
\def\Can{{\mathcal{C}}} 
\def\EL{{\mathcal{EL}}}
\def\Emb{\textnormal{Emb}}
\def\UEmb{\textnormal{UEmb}}
\def\Homeo{\textnormal{Homeo}}
\def\Aff{\textnormal{Aff}}
\def\Q{\mathbb Q}
\def\CP{{\mathbb{CP}}}
\def\F{\mathcal F}
\def\f{\mathfrak f}
\def\X{\mathfrak X}
\def\XX{\mathcal{X}}

\def\length{\textnormal{length}}
\def\area{\textnormal{area}}
\def\rot{\textnormal{rot}}
\def\rott{\textnormal{rot}^\sim}
\def\link{\textnormal{link}}
\def\wind{\textnormal{wind}}
\def\cl{\textnormal{cl}}
\def\scl{\textnormal{scl}}

\newcommand{\marginal}[1]{\marginpar{\tiny #1}}

\title{Sausages}
\author{Danny Calegari}
\address{Department of Mathematics \\ University of Chicago \\
Chicago, Illinois, 60637}
\email{dannyc@math.uchicago.edu}
\date{\today}

\begin{abstract}
The {\em shift locus} is the space of normalized 
polynomials in one complex variable for which every critical point is
in the attracting basin of infinity. The method of {\em sausages} gives
a (canonical) decomposition of the shift locus in each degree
into (countably many) codimension 0
submanifolds, each of which is homeomorphic to a complex algebraic variety.
In this paper we explain the method of sausages, and some of its consequences.
\end{abstract}

\maketitle


\section{Sausages}

For each integer $q\ge 2$ the {\em shift locus} $\SS_q$ is the set of degree $q$ polynomials
$f$ in one complex variable of the form
$$f(z):=z^q + a_2 z^{q-2} + a_3 z^{q-3} + \cdots + a_q$$
for which every critical point of $f$ is in the attracting basin of $\infty$.
One can think of $\SS_q$ as a open submanifold of $\C^{q-1}$; understanding its
topology is a fundamental problem in complex dynamics. 
For example, when $q=2$ the complement of $\SS_2$ in $\C$ is the Mandelbrot set.
Knowing that $\SS_2$ is homeomorphic to a cylinder implies the famous theorem
of Douady--Hubbard that the Mandelbrot set is connected.

Although the $\SS_q$ are highly transcendental spaces, 
the method of {\em sausages} (which we explain in this section) shows that
each $\SS_q$ has a canonical decomposition into codimension 0 submanifolds whose
interiors are homeomorphic to certain explicit algebraic varieties. From this
one can deduce a considerable amount about the topology of $\SS_q$, especially
in low degree.

The construction of sausages has several steps, and goes via an
intermediate construction that associates, to each polynomial $f$ in $\SS_q$,
a certain combinatorial object called a {\em dynamical elamination}. 

\subsection{Green's function}\label{subsection:Greens_function}

Let $K$ be a compact subset of $\C$ with connected complement $\Omega_K:=\C-K$.
If $K$ has positive logarithmic capacity (for example, if the Hausdorff
dimension is positive) then there is a canonical {\em Green's function}
$g:\Omega_K \to \R^+$ satisfying
\begin{enumerate}
\item{$g$ is harmonic;}
\item{$g$ extends continuously to $0$ on $K$; and}
\item{$g$ is asymptotic to $\log{|z|}$ near infinity (in the sense that
$g(z)-\log{|z|}$ is harmonic near infinity).}
\end{enumerate}

There is a unique germ near infinity of a holomorphic function $\phi$, tangent
to the identity at $\infty$, for which $g=\log|\phi(z)|$.

\subsection{Filled Julia Set}

Let $f$ be a degree $q$ complex polynomial. After conjugacy by a complex
affine transformation $z \to \alpha z + \beta$ we may assume that
$f$ is {\em normalized}; i.e.\/ of the form
$$f(z):=z^q + a_2z^{q-2} + a_3z^{q-3} + \cdots + a_q$$
The {\em filled Julia set} $K(f)$ is the set of complex numbers
$z$ for which the iterates $f^n(z)$ are (uniformly) bounded. It is a fact that
$K(f)$ is compact,
and its complement $\Omega_f:=\C-K(f)$ is connected. The union
$\hat{\Omega}_f:=\Omega_f \cup \infty$ is the attracting basin of $\infty$.

B\"ottcher's Theorem (see e.g.\/ \cite{Milnor} Thm.~9.1) 
says that $f$ is holomorphically conjugate near infinity
to the map $z \to z^q$. For normalized $f$, the germ of the
conjugating map $\phi$ (i.e.\/ $\phi$ so that $\phi(f(z))=\phi(z)^q$)
is uniquely determined by requiring that $\phi$ is tangent to the
identity at infinity.  The (real-valued)
function $g(z):=\log|\phi(z)|$ is harmonic, and satisfies the functional
equation $g(f(z))=q\cdot g(z)$. We may extend $g$
via this functional equation to all of $\Omega_f$ and observe that $g$ so
defined is the Green's function of $K(f)$.

\subsection{Maximal domain of $\phi^{-1}$}

Let $\overline{\D}\subset \C$ denote the closed unit disk, and 
$\E:=\C -\overline{\D}$ the exterior. We will use logarithmic
coordinates $h=\log(|z|)$ and $\theta =\arg(z)$ on $\E$ and on Riemann
surfaces obtained from $\E$ by cut and paste. Note that $g=h\phi$ where
$g$ and $\phi$ are as in \S~\ref{subsection:Greens_function}.

\begin{figure}[htpb]
\centering
\includegraphics[scale=0.6]{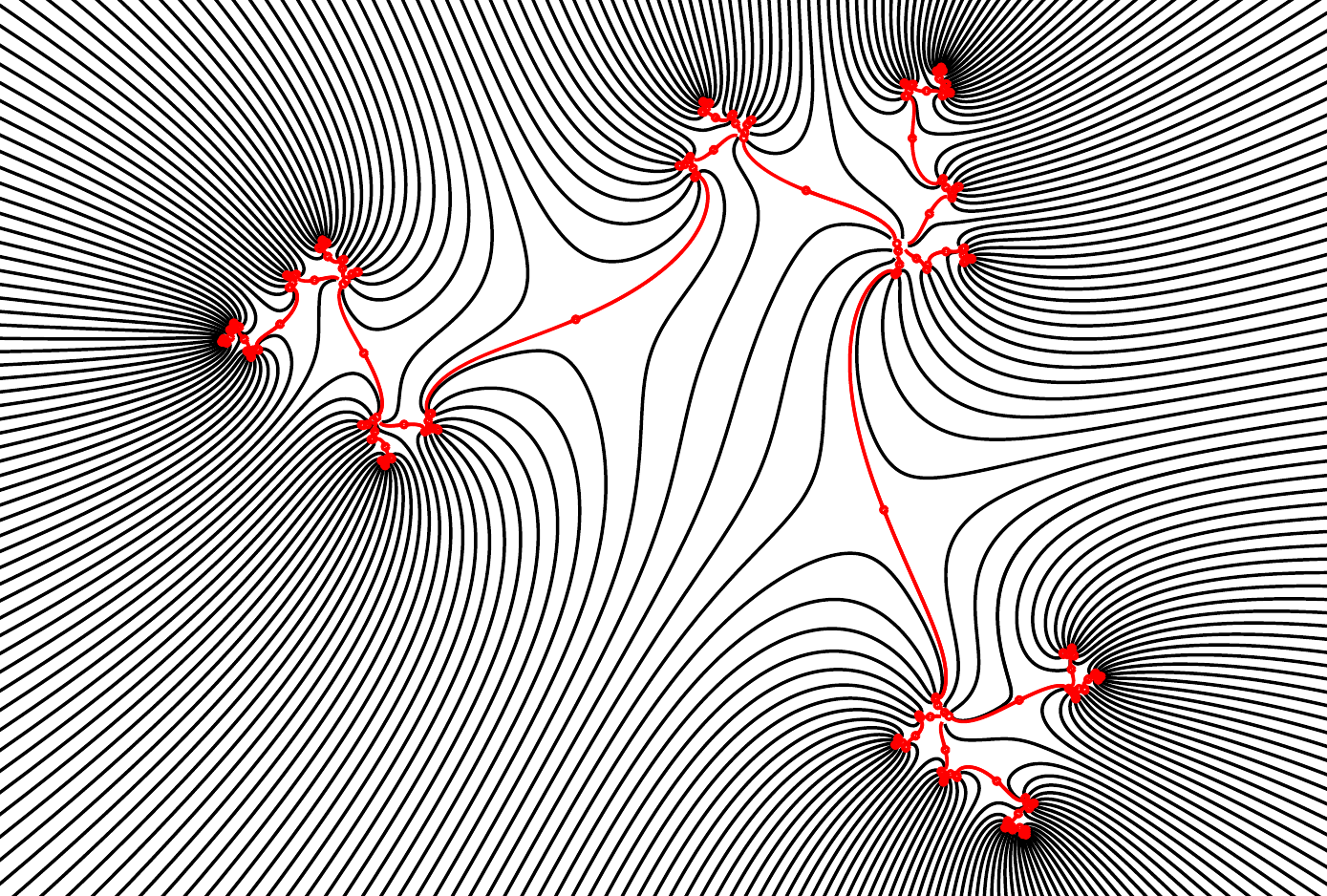} 
\caption{Gradient flowlines of $g$ for a Cantor set $K$}\label{green_flowlines}
\end{figure}

For any $K$ with Green's function $g$ and associated $\phi$ we
can analytically continue $\phi^{-1}$ from infinity along radial lines of $\E$.
The image of these radial lines under $\phi^{-1}$ 
are the descending gradient flowlines of $g$ (i.e.\/ the integral
curves of $-\text{grad}(g)$), and we
can analytically continue $\phi^{-1}$ until the gradient flowlines 
run into critical points of $g$. Figure~\ref{green_flowlines} shows
some gradient flowlines of $g$ for a Cantor set $K$.

Note that some critical points of $g$ might have multiplicity greater
than one; however because $g$ is harmonic, 
the multiplicity of every critical point is finite, and the critical
points of $g$ are isolated and can accumulate (in $\Chat$) only on $K$.
With this proviso about multiplicity, we want to do a sort of `Morse theory'
for the function $g$.

Let $L'$ be the union of the segments of the gradient flowlines of $g$ descending from
all the critical points of $g$; in Figure~\ref{green_flowlines} these 
are in red (gray, for black and white reproduction). 
Then $\Omega_K - L'$ is the image of 
the maximal (radial) analytic extension of $\phi^{-1}$. 
The domain of this maximal extension $\phi^{-1}$ may be described as follows.
For $w \in \E$ define the {\em radial segment} $\sigma(w) \subset \E$ to be the 
set of points $z$ with $\arg(z)=\arg(w)$ and
$|z| \le |w|$. The {\em height} of $\sigma$, denoted $h(\sigma)$, is
$\log(|w|)$. The domain of $\phi^{-1}$ is $\E - L$ where $L$ is
the union of a countable proper (in $\E$) collection of radial segments.

If $K=K(f)$ for a polynomial $f$, the critical points of $g$ are the critical
points {\em and critical preimages} of $f$; i.e.\/ points $z$ for which
$(f^n)'(z)=0$ for some positive $n$.
Thus $L'$ is nearly $f$-invariant: the image
$f(L')$ is equal to $L' \cup \ell'$ where $\ell'$ is the (finite!) set of descending flowlines 
from the critical {\em values} of $f$ in $\Omega_f$
(which are themselves not typically critical).

Likewise the map $z \to z^q$ on $\E$ takes $L$ to $L \cup \ell$ where $\ell$ is a
finite set of radial segments mapped by $\phi^{-1}$ to $\ell'$.

\subsection{Cut and Paste}\label{subsection:cut_and_paste}

\begin{figure}[htpb]
\centering
\includegraphics[scale=0.2]{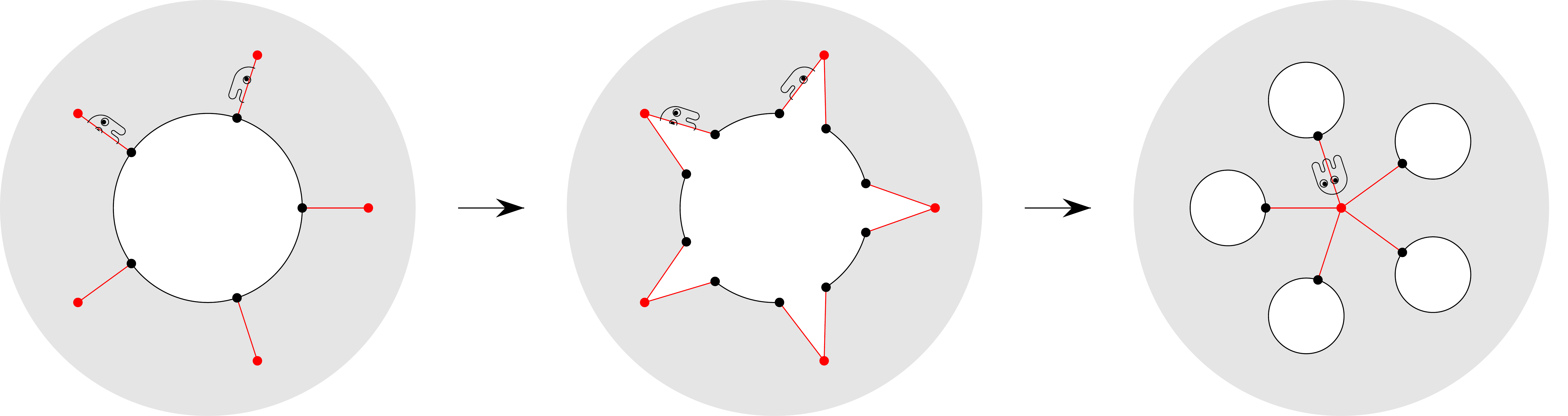} 
\caption{Cut and paste over $L_c$ of multiplicity 4}\label{cut_and_paste}
\end{figure}

Let $c$ be a critical point of $g$ and let $L'_c$ be the union of the
gradient flowlines of $g$ descending from $c$ (and for simplicity, here and in the sequel 
let's suppose these flowlines do not run into another critical point). Then $L'_c$ is the
union of $n+1$ proper embedded rays from $c$ to $K$ where $n$ is the multiplicity
of $c$ as a critical point 
(these rays extend continuously to $K$ when the components of $K$ are
locally connected; otherwise they may `limit to' a {\em prime end} 
of a component of $K$). There is a corresponding collection 
$L_c$ of $n+1$ radial segments
$\sigma_j:=\sigma(w_j)$ all of the same height, 
where indices are circularly ordered according to the 
arguments of the $w_j$. The map $\phi^{-1}$ extends continuously along radial
lines from infinity all the way to the $w_j$: the $w_j$ all map to $c$. 
But any `extension' of $\phi^{-1}$ 
over $L_c$ will be multivalued. We can repair this multivaluedness
by {\em cut and paste}: cut open $\E$ along the segments $L_c$ to create
two copies $\sigma_j^+$ resp. $\sigma_j^-$ for each $\sigma_j$ 
on the `left' resp. `right' in the circular order. 
Then glue each segment $\sigma_j^-$ to $\sigma_{j+1}^+$ by a
homeomorphism respecting absolute value. Under this operation the
collection of segments $L_c$ are reassembled into an `asterisk' which resembles
the cone on $n+1$ points; see Figure~\ref{cut_and_paste}.

The result is a new Riemann surface for which
the map $\phi^{-1}$ now extends (analytically and single-valuedly) over 
the (cut-open and reglued) image of $L_c$, whose image is exactly $L'_c$.

If we perform this cut and paste operation simultaneously 
for all the different $L_c$ making up $L$, the Riemann surface
$\E$ is reassembled into a new Riemann surface $\Omega$
so that $\phi^{-1}$ extends to an {\em isomorphism}
$\phi^{-1}:\Omega \to \Omega_K$.

If $K=K(f)$ for a polynomial $f$, then the map $z \to z^q$ on $\E$ descends to a 
well-defined degree $q$ holomorphic self-map $F:\Omega \to \Omega$ 
and $\phi^{-1}$ conjugates $F|\Omega$ to $f|\Omega_f$.

\subsection{Elaminations}\label{subsection:elaminations}

It is useful to keep track of the partition of $L'$ and $L$ into finite
collections $L'_c$ and $L_c$ associated to the critical points $c$ of $g$.

For each critical point $c$ of multiplicity $n$ we span the $n+1$ segments
of $L_c$ by an ideal hyperbolic $(n+1)$-gon in $\overline{\D}$. The segments
of $L_c$ become the {\em tips} and the ideal polygon becomes the {\em vein}
of a {\em leaf} of multiplicity $n$ in an object called an {\em extended lamination} --- or
{\em elamination} for short. When every critical point has multiplicity 1, we say
the elamination is {\em simple}. See Figures~\ref{dynamical_lamination} and
\ref{tautological} for examples of simple elaminations. 
The key topological property of elaminations
is that the veins associated to different leaves {\em do not cross}. This is
equivalent to the fact that the result $\Omega$ of cut and paste along $L$
is a {\em planar} surface (because it is isomorphic to $\Omega_X \subset \C$).

Elaminations are introduced and studied in \cite{Calegari_sausages}. 
The set $\EL$ of elaminations becomes a space with respect to a certain 
topology (the {\em collision topology}), and can be given the structure of 
a disjoint union of (countable dimensional) complex manifolds. For example,
the space of elaminations with $n-1$ leaves (counted with multiplicity) is 
homeomorphic (but {\em not} biholomorphic) to the space of degree $n$
normalized polynomials with no multiple roots.

\subsection{Dynamical Elaminations}\label{subsection:dynamical_elamination}

\begin{figure}[htpb]
\centering
\includegraphics[scale=0.6]{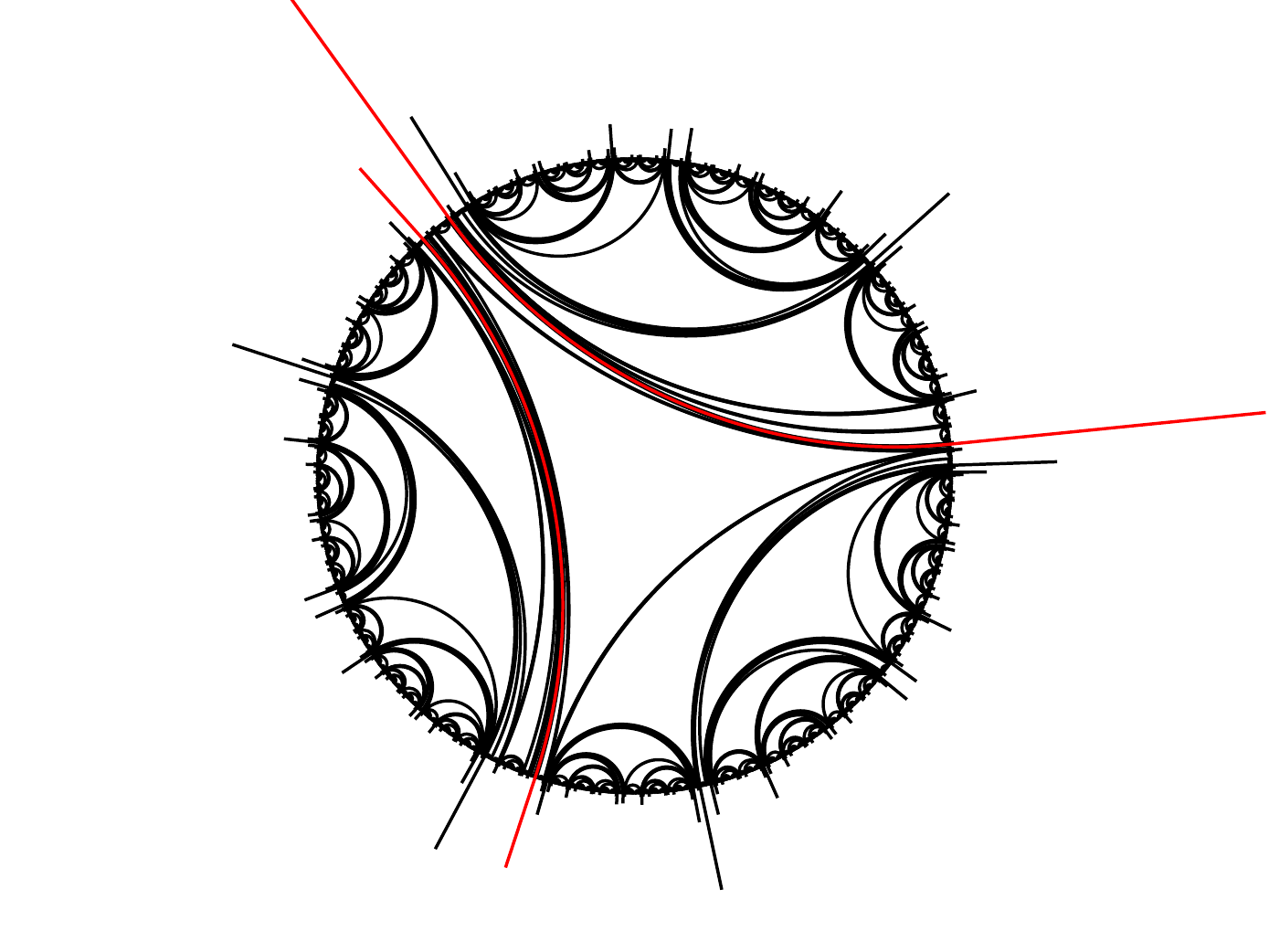} 
\caption{Simple dynamical elamination of degree $3$; critical leaves are in red}\label{dynamical_lamination}
\end{figure}

Figure~\ref{dynamical_lamination} depicts the elamination associated to 
$K(f)$ for a degree 3 polynomial $f$. 
The {\em critical leaves} --- i.e.\/ the leaves with
tips $L_c$ associated to $c$ a critical point of $f$ --- are in red. 
Every other leaf corresponds to a {\em precritical} point of $f$ (which
are critical points of the Green's function). This
elamination is {\em simple}: every leaf has exactly two tips. 

Let $\Lambda$ denote the elamination associated to $L$. Note that $\Lambda$ 
depends not just on $L$ as a set of segments, but also on their partition into
subsets $L_c$.

The map $z \to z^q$ on $\E$ acts on segments and therefore also on leaves, with the
following exception. If $\ell$ is a leaf whose tips have arguments that 
all differ by integer multiples of $2\pi/q$ then these segments will have the
same image under $z \to z^q$. Since leaves should have at least two tips
(by convention), if $\ell$ is a leaf {\em all} of whose tips have
arguments that differ by integer multiples of $2\pi/q$ then the image of
$\ell$ under $z \to z^q$ is undefined.

Suppose $K=K(f)$ for a degree $q$ polynomial. Let $C$ denote the critical leaves
of $L$ (those associated to critical points of $f$). The map $z \to z^q$
takes leaves to leaves in the obvious sense, and takes $\Lambda - C$ to $\Lambda$.

We say an elamination $\Lambda$ is a {\em degree $q$ dynamical elamination}
if 
\begin{enumerate}
\item{it has finitely many leaves $C$ each of whose arguments differ 
by integer multiples of $2\pi/q$ (the {\em critical leaves});}
\item{the map $z \to z^q$ takes $\Lambda - C$ to $\Lambda$; and}
\item{every leaf has exactly $q$ preimages.}
\end{enumerate} 
A degree $q$ dynamical elamination is {\em maximal} if
there are $q-1$ critical leaves, counted with multiplicity.

The elamination $\Lambda$ associated to a degree $q$ polynomial $f$ is a
degree $q$ dynamical elamination. It is maximal if and only if all the
critical points of $f$ are in $\Omega_f$.

A set of (non-crossing) leaves $C$, each with arguments that differ by 
integer multiples of $2\pi/q$ is called a {\em degree $q$ critical set}. A critical
set is {\em maximal} if there are $q-1$ leaves counted with multiplicity.
It turns out that every maximal degree $q$ critical set $C$ is exactly the
set of critical leaves of a {\em unique} (maximal) degree $q$ dynamical elamination
$\Lambda$; see \cite{Calegari_sausages} Prop.~5.3. 
The set of maximal degree $q$ dynamical elaminations is denoted $\DL_q$.
As a subset of $\EL$ it has the structure of an open complex manifold of dimension 
$q-1$ with local coordinates coming from the (endpoints of) segments of $C$ 
(at least at a generic $\Lambda$).

\subsection{The Shift Locus}\label{subsection:shift}

For each degree $q$ the {\em Shift Locus} $\SS_q$ is the space of
degree $q$ normalized polynomials $f(z):=z^q + a_2z^{q-2} + a_3z^{q-3} + \cdots + a_q$
for which every critical point is in the basin of infinity $\Omega_f$.
The coefficients $a_2,\cdots,a_q$ are coordinates on $\SS_q$
realizing it as an open subset of $\C^{q-1}$. 

A polynomial $f$ is in $\SS_q$ if and only if 
the Julia set of $f$ is a Cantor set on which $f$ is uniformly expanding 
(for some metric). Thus for such polynomials, $J(f)=K(f)$ and $\Omega_f$ is
the entire Fatou set (i.e.\/ the maximal domain of normality of $f$ and its 
iterates; see e.g.\/ \cite{Milnor}).

Suppose $f \in \SS_q$ with associated dynamical elamination $\Lambda$. Since
all critical points of $f$ are in $\Omega_f$, it follows that $\Lambda$ is
maximal; thus there is a map $\SS_q \to \DL_q$ called the {\em B\"ottcher map}. 
Conversely, if $\Lambda$ is
any maximal degree $q$ dynamical elamination, and $\Omega$
is obtained from $\E$ by cut and paste along $\Lambda$, then $F|\Omega$
extends (topologically) over the space of ends of $\Omega$ to define a
degree $q$ self-map $\bar{F}$ of a topological sphere $\bar{\Omega}\cong S^2$. It turns
out that there is a canonical conformal structure on $\bar{\Omega}$ extending
that on $\Omega$ so that $\bar{F}$ is holomorphic. After choosing 
suitable coordinates on $\bar{\Omega}$ near $\infty$ the map $\bar{F}$ 
becomes a degree $q$ normalized polynomial, which is contained in $\SS_q$. 
The analytic content of this theorem is essentially due to de Marco--McMullen; 
see e.g.\/ \cite{de_Marco_McMullen} Thm.~7.1 or 
\cite{Calegari_sausages} Thm.~5.4 for a different proof.

Thus the B\"ottcher map $\SS_q \to \DL_q$ is a homeomorphism 
(and in fact an isomorphism of complex manifolds). 

\subsection{Stretching and Spinning}\label{subsection:stretch}

There is a (multiplicative) $\R^+$ action on $\EL$ called {\em stretching}
where $t \in \R^+$
acts on $\Lambda$ by multiplying the $h$ coordinate of every leaf by $t$. This
action is free and proper. It preserves $\DL_q$ for each $q$, and shows 
that $\DL_q$ (and therefore also $\SS_q$) is homeomorphic to the
product of $\R$ with a manifold of (real) dimension $2q-3$.
It is convenient for what follows to define $\DL_q'$ to be the open subspace
of $\DL_q$ for which the highest critical leaf has $\log_q(h) \in (-1/2,1/2)$.
By suitably `compressing' orbits of the $\R^+$ action we see there is a
homeomorphism $\DL_q \to \DL_q'$.

There is also an $\R$ action on $\EL$ called {\em spinning} 
where $t \in \R$ simultaneously 
rotates the arguments of leaves of height $h$ by $th$. 
This makes literal sense for the
(finitely many) leaves of greatest height. When leaves of lesser height are
collided by those of greater height the shorter leaf is `pushed over' the taller
one; the precise details are explained in \cite{Calegari_sausages} \S~3.2.
This $\R$ action also preserves each $\DL_q$. The closure
of the $\R$-orbits in each $\DL_q$ are real tori, 
and the $\R$-orbits sit in these tori
as parallel lines of constant slope. A typical 
$\R$-orbit has closure which is a torus of real dimension $q-1$, but if
some critical leaves have multiplicity $>1$ or if
distinct critical leaves have rationally related heights, the closure will
be a torus of lower dimension. 

Stretching and spinning combine to give an action of
the (oriented) affine group $\R \rtimes \R^+$ of the line on $\EL$ and on
each individual $\DL_q$.

\subsection{Sausages}

Suppose $K=K(f)$ for a degree $q$ polynomial. The map $f$ is algebraic, but the
domain $\Omega_f$ is transcendental. When we move to the elamination side,
the map $z\to z^q$ and the domain $\E$ are (semi)-algebraic, but the
combinatorics of $L$ is hard to understand. Sausages are a way to find a substitute
for $(f,\Omega_f)$ for which both the map and the domain are algebraic and more 
comprehensible.

The idea of sausages is to find a dynamically invariant way to cut up
the domain $\Omega$ into a {\em tree of Riemann spheres}, so that $F$ induces
polynomial maps between these spheres. The sausage map is {\em not} holomorphic,
but it induces homeomorphisms between certain codimension 0 submanifolds of
$\DL_q'$ and certain explicit algebraic varieties whose topology is in some
ways much easier to understand.

Now let's discuss the details of the construction.
First consider the map $z\to z^q$ on $\E$ alone. Let $h:=\log(|z|)$
and $\theta=\arg(z)$ be cylindrical coordinates on $\E$, so that $\E$
becomes the half-open cylinder $S^1 \times \R^+$ in $(\theta,h)$-coordinates,
and $z\to z^q$ becomes the map which is multiplication by $q$ which we denote $\times q$. 
For each integer $n$ let $I_n$ denote the open interval $(q^{n-1/2},q^{n+1/2})$ and
let $A_n$ be the annulus in $\E$ where $h \in I_n$ and let 
$A=\cup_n A_n \subset \E$; the complement of $A$
is the countable set of circles with $\log_q(h) \in 1/2 + \Z$.
Then $\times q$ takes $A_n$ to $A_{n+1}$. 

This data is holomorphic but not algebraic. So let's choose (rather
arbitrarily) an orientation-preserving diffeomorphism $\nu_0:I_0 \to \R$
and for each $n$ define $\nu_n:I_n \to \R$ by 
$\nu_n(h) = q^n\nu_0(q^{-n}h)$ (so that by induction the $\nu_n$ satisfy
$\nu_{n+1}(qh)=q\nu_n(h)$ for all $n$ and $h \in I_n$), 
and define $\mu:A \to S^1 \times \R$ to be
the map that sends $(\theta,h)$ to $(\theta,\nu_n(h))$ if $(\theta,h) \in A_n$.
By construction, $\mu$ commutes with multiplication by $q$:
$$\mu(q\theta,qh) = (q\theta,\nu_{n+1}(qh)) = (q\theta,q\nu_n(h))=q\mu(\theta,h)$$
In other words, $\mu$ semi-conjugates $\times q$ on $A$ to $\times q$ on
$S^1 \times \R$, which (by exponentiating) becomes the map $z \to z^q$ on $\C^*$, an algebraic
map on an algebraic domain. Actually, it is better to keep a separate copy
$\C^*_n:=\mu(A_n)$ of $\C^*$ for each $n$, so that $\mu$ {\em conjugates}
$\times q$ on $A$ to the self-map of $\cup_n \C^*_n$
which sends each $\C^*_n$ to $\C^*_{n+1}$ by $z \to z^q$.

\subsection{Sausages and Dynamics}

\begin{figure}[htpb]
\centering
\includegraphics[scale=0.08]{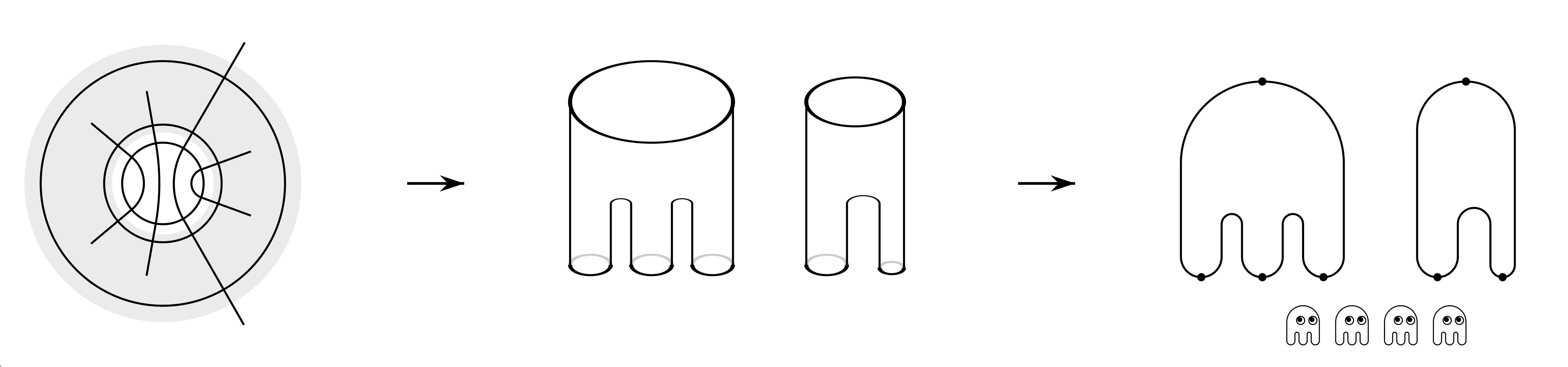} 
\caption{$A_n$ is cut and paste into $B_n$ which in turn maps to a disjoint union of Riemann spheres}\label{sausage_example}
\end{figure}

Now suppose we have a dynamical elamination $\Lambda$ with critical leaves $C$ invariant
under $z \to z^q$. For each $A_n$ the tips of $\Lambda$ intersect $A_n$ 
in a finite collection of vertical
segments $L_n$ (some of which will pass all the way through $A_n$) and we can
perform cut-and-paste separately on each $A_n$ to produce a (typically disconnected)
surface $B_n$. Furthermore
we can perform cut-and-paste on $\C^*_n$ along the image $\mu(L_n)$ which,
by construction, is compatible with the Riemann surface structure. The
result is to cut and paste $\C^*_n$ into a finite collection of algebraic 
Riemann surfaces, each individually isomorphic to $\C$ minus a 
finite set of points and which can be canonically completed to 
Riemann spheres in such a way that the map $F$ on $\Omega$ descends to
a map $\f$ from this union of Riemann spheres to itself; see Figure~\ref{sausage_example}.

Denote the individual Riemann spheres by $X_v$, and by abuse of notation,
write $\f_v:X_v \to X_{\f(v)}$ for the restriction of $\f$ to the 
component $X_v$. By the previous discussion, each map $\f_v$ is {\em holomorphic},
so that if we choose suitable coordinates on $X_v$ and $X_{\f(v)}$ the map
$\f_v$ becomes a polynomial. There is almost a canonical choice of
coordinates, which we explain in the next two sections.

Each $X_v$ corresponds to a component $B_v$ of some $B_n$, and
gets a canonical finite set of {\em marked points} $Z_v'$ which correspond to
the `boundary circles' of $B_v$. The unique boundary
circle with greatest $h$ coordinate picks out a point that we can identify
with $\infty\in X_v$; we denote by $Z_v$ the set consisting of the rest of the
marked points. 
The collection of individual Riemann spheres $X_v$
can be glued up along their marked points into an infinite 
genus zero nodal Riemann surface so that the indices $v$ are
parameterized by the vertices $v$ of the tree of gluings $T$. This tree
is oriented, so that an edge $v$ goes to $w$ if $X_v$ is glued along
$\infty$ to one of the (non-infinite) marked points of $X_w$. We call
$w$ the {\em parent} of $v$ and $v$ one of the {\em children} of $w$. 
If we make the assumption that no boundary component of any $B_v$ contains
a critical point (this is the generic case) then each $\zeta \in Z_w \subset X_w$ is
attached to a unique $X_v$ for $v$ some child of $w$.
If $v$ is a child of $w$, and $X_v$ is glued to $X_w$ at the point 
$\zeta \in Z_w$, then if $\zeta$ is a critical point of $\f_w$ of multiplicity
$m$, the degree of $\f_v$ is $m+1$. By abuse of notation, we denote 
the induced (simplicial, orientation-preserving) map on $T$ also by $\f$.

If $\Lambda$ is empty, then $T$ is just a line, and each vertex has
a unique child. If $\Lambda$ is nonempty, then since there are only finitely
many leaves of greatest height,
there is a unique highest vertex $v$ of $T$
with more than one child. Let $w$ be the parent of $v$.
The uppermost boundary components of $B_v$ and $B_w$ are
canonically identified with the unit circle $S^1:=\R/2\pi\Z$. 
By identifying these circles
with the unit tangent circles at $\infty$ in $X_v$ and $X_w$ we can choose
coordinates on these Riemann spheres so that the tangent to the positive
real axis corresponds to the angle $0 \in S^1$. In these coordinates 
$X_v$ and $X_w$ are identified with copies $\Chat_v$ and $\Chat_w$ of the Riemann
sphere $\Chat$, 
and after precomposing with a suitable complex affine translation, 
$\f_v$ becomes a normalized degree $q$ polynomial map 
$\f_v:z \to z^q + b_2 z^{q-2} + \cdots b_q$, and the (non-infinite) marked
points of $X_v$ become the roots of $\f_v$ in $\Chat_v$.

Vertices of $T$ above $v$ and the maps between their respective Riemann
surfaces do not carry any information. Let $w_1:=w$ denote the parent of
$v$, and inductively let $w_n$ be the parent of $w_{n-1}$. Then each
$X_{w_n}$ has exactly two marked points, which we can canonically identify
with $\infty$ and $0$, and the map $\f_{w_{n-1}}:\Chat_{w_{n-1}} \to \Chat_{w_n}$
is canonically normalized as $z \to z^q$. 

Since these vertices carry no information, we discard them.
Thus we make the convention that $T$ is the {\em rooted} tree consisting of
$v$ together with its (iterated) children, and we let $\X$ denote the 
nodal Riemann surface corresponding
to the union of $X_w$ with $w$ in $T$. We record the data of the polynomial
$\f_v$ associated to the root $v$, 
though we do not interpret this any more as a map between Riemann spheres, so
that $\f$ is now a map from $\X-X_v$ to $\X$ and $\f_v$ is a polynomial
function on $X_v \cong \Chat$.

\subsection{Tags and sausage polynomials}\label{subsection:tags}

The choice of a distinguished point on a boundary $S^1$ component of some 
$B_u$ is called a {\em tag}. Tags are the data we need to choose coordinates on
$\X$ so that every $\f_u$ becomes a polynomial.
We may identify this boundary circle with the unit tangent circle at a
marked point on $X_u$, and think of the tag as data on $X_u$. By induction, we can
choose tags on $X_u$ in the preimage of the tags of $X_{\f(u)}$ under the
map $\f_u:X_u \to X_{\f(u)}$ and inductively define 
coordinates $\Chat_u$ on $X_u$ for which $\f_u$ is represented by a 
normalized polynomial map (in general of degree $\le q$).

Suppose $u$ has parent $u'$, and $\infty$ in $\Chat_u$ is attached at some
point $\zeta \in Z_{u'} \in \Chat_{u'}$. Suppose $\zeta$ is a critical point of $\f_{u'}$ 
with multiplicity $m$. Then $\f_u$ has degree $m+1$. There are $m+1$
different choices of tag at $\zeta$ that map to the tag at
$\f_{u'}(\zeta)$, and the different choices affect the normalization of
$\f_u$ by precomposing with multiplication by an $(m+1)$st root of unity.

The endpoint of this discussion is that we can recover $\X,\f$ from the data of
a rooted tree $T$, and a set of equivalence classes of pair 
$(\text{tag},\text{normalized polynomial }\f_u)$. Call this data a
(degree $q$) {\em sausage polynomial}.

A dynamical elamination $\Lambda$ is {\em generic} if the critical points of $F$ are all
contained in $A$; i.e.\/ if no critical (or by induction, precritical) point
has $h$ coordinate with $\log(h) \in 1/2 + \Z$. The {\em sausage map} is the map 
that associates a sausage polynomial to a degree $q$ dynamical elamination. 
A sausage polynomial is {\em generic} 
resp. {\em maximal} if it is in the image of a generic resp. maximal dynamical elamination. 

A polynomial $\f_w$ associated to a (generic) sausage polynomial has two kinds
of critical points. The {\em genuine} critical points are those in
$\Chat_w - Z_w'$ (recall that $Z_w'$ is $Z_w \cup \infty$). 
The {\em fake} critical points are those in $Z_w'$ ($=\infty \cup Z_w$) which
correspond to circle components of $B_w$ mapping with degree $>1$.
For a generic dynamical elamination, the genuine critical points of the
associated sausage polynomial are exactly the images of the critical points of
the elamination (i.e.\/ the endpoints of the critical leaves) under the sausage map. 
Thus for a generic maximal sausage
polynomial of degree $q$, there are exactly $q-1$ genuine critical points,
counted with multiplicity.

For a generic, maximal sausage polynomial,
all but finitely many $\f_v$ have degree one. A degree
one map uniquely pulls back tags, and has only one possible normalized
polynomial representative, namely the identity map $z \to z$. Thus a 
generic, maximal sausage polynomial is described by a finite amount of 
combinatorial data together with a finite collection of normalized polynomials.
The reader who wants to see some examples should look
ahead to \S~\ref{subsection:degree_2} and \S~\ref{subsection:degree_3}.

Let $\XX_q$ denote the space of generic maximal degree $q$ sausage polynomials.
Then $\XX_q$ is the disjoint union of countably infinitely many components,
indexed by the combinatorics of $T$ and the degrees of the
normalized polynomials between the associated Riemann spheres. Each
component of $\XX_q$ is a {\em quasiprojective complex variety}
of complex dimension $q-1$. In fact, each component is an iterated fiber
bundle whose base and fibers are certain {\em affine} (complex) varieties
called {\em Hurwitz varieties}, which we shall describe in more detail in
\S~\ref{subsection:Hurwitz_variety}.

\subsection{Sausage space}\label{subsection:sausage_isomorphism}

Recall that
$\DL_q' \subset \DL_q$ denotes the set of maximal degree $q$ dynamical elaminations 
for which the highest critical point has $\log_q(h) \in (-1/2,1/2)$. Let
$\DL_q'' \subset \DL_q'$ denote the subspace of {\em generic} maximal degree
$q$ dynamical elaminations. Then the construction of the previous few sections
defines a map $\DL_q'' \to \XX_q$.

In fact, this map is invertible. Given a sausage polynomial $\X,\f$ over
a tree $T$ with root $v$ we can inductively construct (singular) vertical
resp. horizontal foliations on each $\Chat_w$ as follows. 
On $\Chat_v$ we pull back the foliations of $\C^*$ by lines resp. circles of
constant argument resp. absolute value under the polynomial $\f_v$. Then
on every other $w$ we inductively pull back these foliations under 
$\f_w:\Chat_w \to \Chat_{\f(w)}$. These foliations all carry coordinates
pulled back from $\C^*$, and $\Chat_w$ minus infinity 
and its marked points becomes isomorphic to a branched Euclidean
Riemann surface with ends isomorphic to the ends of (infinite) 
Euclidean cylinders. We can reparameterize the vertical coordinates
on each of these Riemann surfaces by the inverses of the maps $\nu_n$, and
then glue together the result by matching up boundary circles using tags.
This defines an inverse to the map $\DL_q'' \to \XX_q$ and shows that
this map is a homeomorphism. See \cite{Calegari_sausages} Thm.~9.20 for details.

\subsection{Decomposition of the Shift Locus}

Putting together the various constructions we have discussed so far we
obtain the following summary:

\begin{enumerate}
\item{\S~\ref{subsection:shift}: The map that associates to $f\in \SS_q$ a maximal degree $q$ dynamical
elamination $\Lambda$ gives an isomorphism of complex manifolds
$\SS_q \to \DL_q$.}
\item{\S~\ref{subsection:stretch}: By compressing orbits of the free $\R^+$ action on $\DL_q$ 
we obtain a homeomorphism
$\DL_q \to \DL_q'$ to the subspace whose largest critical leaf has log height
$\log_q(h)\in (-1/2,1/2)$.}
\item{\S~\ref{subsection:sausage_isomorphism}: The open dense subset $\DL_q'' \subset \DL_q'$ 
of generic dynamical elaminations maps
homeomorphically by the sausage map $\DL_q'' \to \XX_q$.}
\item{\S~\ref{subsection:tags}: The space $\XX_q$ is the disjoint union of countably many 
quasiprojective complex varieties, each of which has the structure of an
iterated bundle of affine (Hurwitz) varieties.}
\end{enumerate}

In words: the shift locus $\SS_q$ in degree $q$ has a canonical
decomposition into codimension 0 submanifolds whose
interiors are homeomorphic to certain explicit algebraic varieties. It is
a fact that we do not explain here (see \cite{Calegari_sausages} \S~8
especially Thm.~8.11) that the
abstract cell complex which combinatorially parameterizes the decomposition of 
$\SS_q$ into these pieces is {\em contractible}, so that all the interesting
topology of $\SS_q$ is localized in the components of $\XX_q$.

In the remainder of the paper we give examples, and explore some of
the consequences of this structure.

\section{Sausage Moduli}

Each component $Y$ of $\XX_q$ parameterizes sausages of a fixed combinatorial type.
The combinatorial type determines finitely many vertices $u$ for which the
(normalized) polynomial $\f_u$ has degree $>1$. The combinatorics
constrains these polynomials by imposing conditions on their critical values,
for instance that the critical values are required to lie outside a certain
(finite) set. Thus, each component has the structure of an algebraic variety
which is an iterated fiber bundle, and so that the base and each fiber is
something called a {\em Hurwitz Variety}.

For this and other reasons, the spaces $\SS_q$ and the components 
$Y$ of which they are built bear a close family resemblance to  
the kinds of discriminant complements that arise in the study of
classical braid groups. The full extent of this resemblance is an open
question, partially summarized in Table~\ref{finitetable}.

\subsection{Degree 2}\label{subsection:degree_2}

\begin{figure}[htpb]
\centering
\includegraphics[scale=0.3]{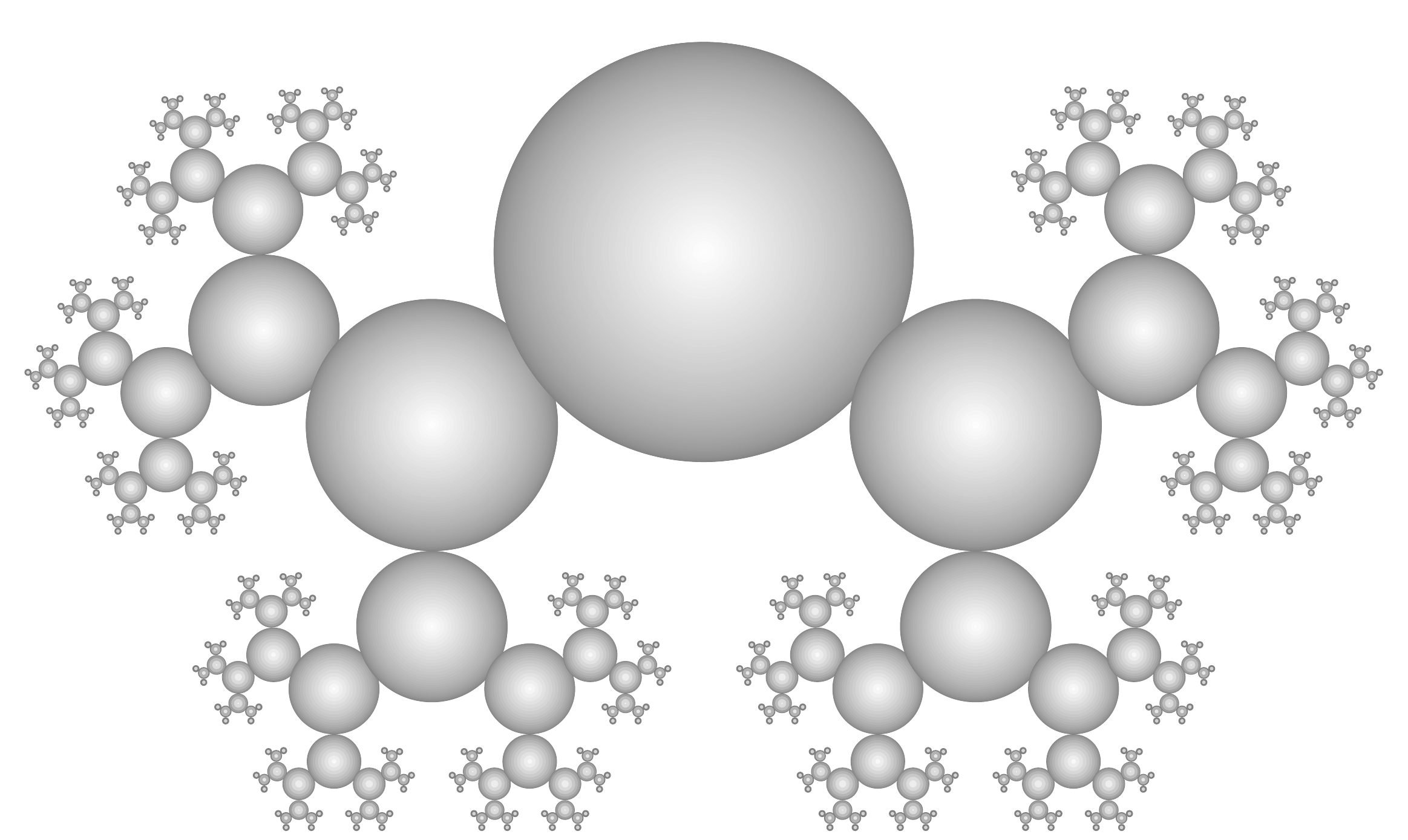} 
\caption{A degree 2 sausage; each vertex is attached to its parent at the
points $\pm \sqrt{-c}$.}\label{degree_2_sausage}
\end{figure}

Let $\X,\f$ be a generic maximal sausage polynomial of degree $2$. 
The root polynomial $\f_v$ is quadratic and normalized.
It has one critical point, necessarily genuine. Thus $\f_v(z):= z^2 + c$
for some $c \ne 0$. Every other vertex $w$ has a polynomial $\f_w$ of degree one;
since polynomials are normalized, $\f_w(z):=z$. Thus all the information
is contained in the choice of the (nonzero) constant coefficient $c$ of
$\f_v$, so that $\XX_2=\C^*$.
The tree $T$ is an infinite dyadic rooted tree, where every
vertex is attached to its parent at the points $\pm \sqrt{-c}$; 
see Figure~\ref{degree_2_sausage}. 

Furthermore, in this case $\DL_2'=\DL_2''$ so that $\SS_2$ is homeomorphic
(but {\em not} holomorphically isomorphic) to $\C^*$. As a corollary one
deduces the famous theorem of Douady--Hubbard \cite{Douady_Hubbard}
that the Mandelbrot set $\M$ (i.e.\/ $\C-\SS_2$) is connected.

\subsection{Discriminant Locus}\label{subsection:discriminant}

In any degree $q$ there is a unique component of $\XX_q$ for which all the
(genuine) critical points are in the root vertex. Thus $\f_v$ is a degree $q$
normalized polynomial with no fake critical points. Since the marked points
$Z_v$ of the root vertex are exactly the {\em roots} of $\f_v$, this means
that $\f_v$ is a normalized polynomial with no critical roots. Equivalently,
$\f_v$ has $q$ distinct roots, so that $\f_v$ is in $Y_q:=\C^{q-1}-\Delta_q$ where
$\Delta_q$ is the {\em discriminant locus}. As is well-known, $Y_q$
is a $K(B_q,1)$ where $B_q$ denotes the braid group on $q$ strands.

\subsection{Degree 3}\label{subsection:degree_3}

Let $\X,\f$ be a generic maximal sausage polynomial of degree $3$. If the root
polynomial $\f_v$ has two genuine critical points we are in the case
discussed in \S~\ref{subsection:discriminant} and the corresponding
component of $\XX_3$ is a $K(B_3,1)$. Otherwise, since the 
root polynomial must have at least one genuine critical point, if it does not have
two it must have exactly one and $\f_v$ is of the form 
$z \to (z-c)^2(z+2c)$ for some $c \in \C^*$.

\begin{figure}[htpb]
\centering
\includegraphics[scale=0.5]{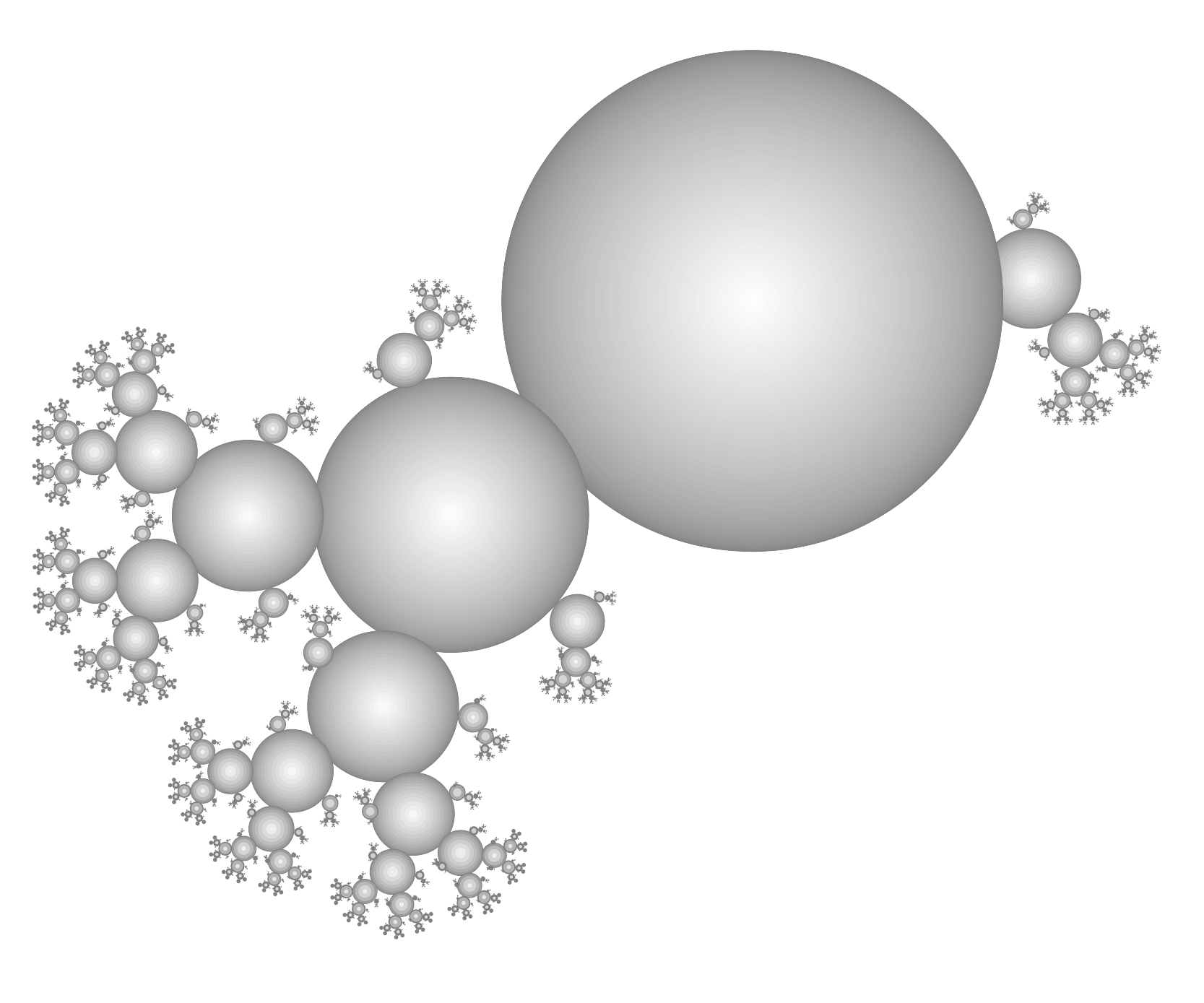} 
\caption{A degree 3 sausage; the root $v$ has $Z_v:= \lbrace c, -2c \rbrace$. 
The child $w_1$ has $Z_{w_1}:=\lbrace \pm \sqrt{c-d}, \pm \sqrt{-2c-d}\rbrace$}\label{degree_3_sausage}
\end{figure}

The (finite) marked points $Z_v$ of $\Chat_v$ are $c$ and $-2c$, and the root
vertex correspondingly has two children $w_1,w_2$ where $\Chat_{w_1}$ is 
attached at $c$ and $\Chat_{w_2}$ is attached at $-2c$. Because $c$ is a double
root, the polynomial $\f_{w_1}$ has degree $2$; because $-2c$ is a simple root,
the polynomial $\f_{w_2}$ has degree $1$.

Write $\f_{w_1}:z \to z^2 + d$. If $d \ne c,-2c$ then $Z_{w_1}$ has four
(non-critical) points (the distinct square roots of $c-d$ and $-2c-d$) and
every other $\f_u$ is degree 1. See Figure~\ref{degree_3_sausage}. 
Thus $c$ and $d$ are moduli parameterizing
a single component of $\XX_3$, and topologically this component is 
a bundle over $\C^*$ whose fiber is homeomorphic to $\C-\lbrace c,-2c\rbrace$.

If $d=c$ or $d=-2c$ then $0$ is a fake critical point for $\f_{w_1}$, and if
$u$ is the child of $w_1$ for which $\Chat_u$ is attached at $0$ 
then $\f_u$ has degree $2$. Since $\f$ is maximal, there is always some vertex
$u'$ at finite combinatorial distance from the root for which $\f_{u'}$ has degree
$2$ and for which the critical point $0$ of $\f_{u'}$ is genuine. Thus each
component of $\XX_3$ is a bundle over $\C^*$ with fiber homeomorphic to $\C$
minus finitely many points.

\subsection{The Tautological Elamination}

The combinatorics of the components of $\XX_3$ is quite complicated. Each component of
$\XX_3$ (other than the discriminant complement c.f.\/ \S~\ref{subsection:discriminant}) is a punctured
plane bundle over the curve $\C^*$ with parameter $c$, and these components
glue together in $\SS_3$ to form a bundle over $\C^*$ whose fiber $\Omega_T$ is 
homeomorphic to a plane minus a Cantor set. 

Actually, there is another
description of $\Omega_T$ in terms of elaminations. For each degree 3 critical 
leaf $C$ there is a certain elamination $\Lambda_T(C)$ called the 
{\em tautological elamination} which can be defined as follows. Let's suppose
that we have a maximal degree $3$ dynamical elamination with two critical
leaves $C$ and $C'$, and that $C$ has the greater height. If we fix $C$, then
$\Omega_T$ parameterizes the space of configurations of $C'$.

The elamination $\Lambda_T(C)$ is defined as follows. With $C$ fixed, each 
choice of (noncrossing) $C'$ determines a dynamical elamination $\Lambda$.
By hypothesis $h(C')< h(C)$ and there are only finitely many (perhaps zero)
precritical leaves $P$ of $C$ with $h(P)>h(C')$. As we vary $C'$ the laminations
$\Lambda$ also vary (in rather a complicated way) but while $h(P)>h(C')$ the
leaves $P$ stay fixed under continuous variations of $C'$. It might happen
that as we vary the leaf $C'$ it collides with a leaf $P$ with $h(P)>h(C')$; the 
elamination $\Lambda_T(C)$ consists of the {\em cubes} $P^3$ of all such $P$
(there is a similar, though more complicated construction in higher degrees).
The fact that $\Lambda_T(C)$ is an elamination is not obvious from this definition.

\begin{figure}[htpb]
\centering
\includegraphics[scale=0.5]{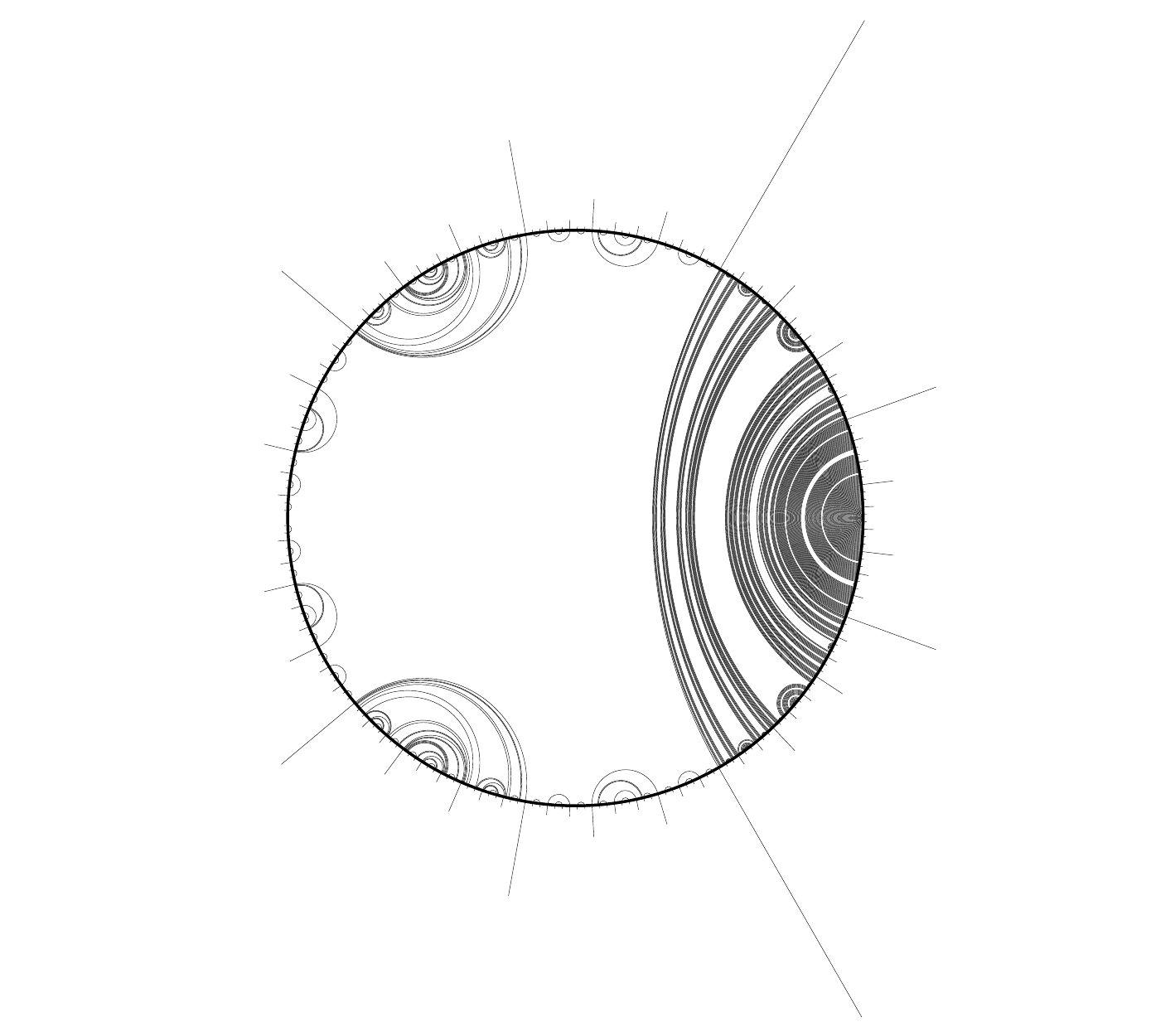} 
\caption{The tautological elamination $\Lambda_T(C)$ for 
$\arg(C)=\pm\pi/6$}\label{tautological}
\end{figure}

The result of cut and paste (as in \S~\ref{subsection:cut_and_paste}) 
on the annulus $1<|z|<|C|$
(thought of as a subset of $\E$) along $\Lambda_T(C)$ is a Riemann surface
$\Omega_T(C)$ holomorphically isomorphic to the moduli space of degree
3 maximal dynamical elaminations for which $C$ is the unique critical leaf of
greatest height. Figure~\ref{tautological} depicts the elamination $\Lambda_T(C)$
for a particular value of $C$ whose tips have angles $\pm \pi/3$.
 
These $\Omega_T(C)$ are the leaves of a (singular) 
one complex dimensional holomorphic foliation of $\SS_3$.

Although it is not a dynamical elamination, the tautological elamination
$\Lambda_T(C)$ is in a natural way the increasing union of 
finite elaminations $\Lambda_n$, namely the leaves of the form $P^3$ as above
where $P$ is a depth $n$ preimage of $C$. Let $\overline{\E}$ denote the closure of
$\E$ in $\C$ so that $\overline{\E}=\E\cup S^1$, the union of $\E$ with the
unit circle. The result $\Omega_n$ of cut and pasting $\E$ along $\Lambda_n$
is partially compactified by a finite set of circles, obtained from $S^1$.
By abuse of notation we denote this finite set of circles by 
$S^1 \mod \Lambda_n$. It turns out that the 
the components of $\XX_3 \cap \Omega_T$ corresponding to sausage polynomials
with fixed $c \in \C^*$ and 
for which the second genuine critical point is in a vertex at depth $n+1$ are
in bijection with the set of components of $S^1 \mod \Lambda_n$.
In fact, more is true.

For each combinatorial type $\X,\f$ let $u$ be the vertex containing the
second genuine critical point (the first, by hypothesis, is contained in the root). 
We define the {\em depth} $n$ of $\X,\f$ to be
the combinatorial distance of $u$ to the root. There is another 
invariant of $\X,\f$: the {\em $\ell$-value}, defined as follows. 
Under iteration of $\f$ (acting on the tree) the vertex $u$ has a length $n$ orbit 
terminating in the root (note that $\f(u)$ is {\em not} typically equal 
to the parent of $u$, but it does have the same depth as the parent). 
The point $\infty$ in $\Chat_u$ is mapped to $\infty$ in $\Chat_{\f(u)}$ and
so on. The product of the degrees of the
polynomials $\f_{\f^i(u)}$ up to but not including the root is some power
of $2$; by definition, $\ell$ is this number divided by $2$. 
The invariants $n$ and $\ell$,
taking discrete values, are really invariants of the components of 
$\XX_3$ and ipso facto of the components of $\XX_3 \cap \Omega_T$.

Here is the relation to $\Lambda_T(C)$. Components of $\XX_3 \cap \Omega_T$ of depth 
$n+1$ are in bijective correspondence with components of $S^1 \mod \Lambda_n$, and
a component of $\XX_3 \cap \Omega_T$ with $\ell$-value $\ell$ corresponds
to a component of $S^1 \mod \Lambda_n$ of length $2\pi \ell \cdot 3^{-n}$.

\subsection{Combinatorics}

Let $N_3(n,m)$ denote the number of components of $S^1$ mod $\Lambda_n$ with depth
$n+1$ and $\ell = 2^m$. We do not know a simple closed form for $N_3(n,m)$ and perhaps
none exists --- one subtle issue is that there are several combinatorially
different ways that a component can have a particular $\ell$-value. However,
an $\ell$-value of $1$ is special, since it corresponds to an $\f$ for which
$\f_{\f^i(u)}$ has degree $1$ for all positive $i$. Correspondingly, there is
an explicit formula for $N_3(n,0)$ that we now give; see
\cite{Calegari_combinatorics} Thm.~3.6 for a proof.

First of all, $N_3(n,0)$ satisfies the recursion
$N_3(0,0)=1$, $N_3(1,0)=1$ and 
$$N_3(2n,0)=3N_3(2n-1,0) \text{ and } N_3(2n+1,0)=3N_3(2n,0)-2N_3(n,0)$$
Knowing this, one can write down an explicit generating function for 
$N_3(n,0)$; the generating function is $(\beta(t)-1)/3t$ where 
$$\beta(t)=\Bigl(\sum_{n=0}^\infty h(n)t^n\Bigr) \prod_{j=0}^\infty \frac {1}{(1-3t^{2^j})}$$
and where the numbers $h(n)$ are defined by
$$h(0)=1 \text{ and } h(n) = (-3)^{s(n)}(1- (-2)^{k(n)})$$
where $2^{k(n)}$ is the biggest power of 2 dividing $n$, and $s(n)$ is the
sum of the binary digits of $n$.

Table~\ref{nltable} gives values of $N_3(n,m)$ for $0 \le n,m \le 12$. Note that $N_3(n,m)=0$
for $n/2 < m < n$; see \cite{Calegari_combinatorics} Thm~5.9.

\begin{table}[ht]
{\tiny
\centering
\begin{tabular}{c|c c c c c c c c c c c c c} 
 $n\backslash \ell$ & 1 & 2 & $2^2$ & $2^3$ & $2^4$ & $2^5$ & $2^6$ & $2^7$ & $2^8$ & $2^9$ & $2^{10}$ & $2^{11}$ & $2^{12}$\\
 \hline
 0 &  1\\ 
 1 &  1 & 1\\
 2 &  3 & 1 & 1\\
 3 &  7 & 6 & 0 & 1\\
 4 &  21 & 16 & 3 & 0 & 1\\
 5 &  57 & 51 & 13 & 0 & 0 & 1\\
 6 &  171 & 149 & 39 & 5 & 0 & 0 & 1\\
 7 &  499 & 454 & 117 & 23 & 0 & 0 & 0 & 1\\
 8 &  1497 & 1348 & 360 & 66 & 9 & 0 & 0 & 0 & 1\\
 9 &  4449 & 4083 & 1061 & 207 & 41 & 0 & 0 & 0 & 0 & 1\\
 10 & 13347 & 12191 & 3252 & 591 & 126 & 17 & 0 & 0 & 0 & 0 & 1\\
 11 & 39927 & 36658 & 9738 & 1799 & 370 & 81 & 0 & 0 & 0 & 0 & 0 & 1\\
 12 & 119781 & 109898 & 29292 & 5351 & 1125 & 240 & 33 & 0 & 0 & 0 & 0 & 0 & 1\\
\end{tabular}
}
\vskip 12pt
\caption{Number of components of length $\ell/3^n$ at depth $n$}\label{nltable}
\end{table}

\subsection{Hurwitz Varieties}\label{subsection:Hurwitz_variety}

Let $X$ be a component of $\XX_q$ parameterizing sausage polynomials of a fixed
combinatorial type. $X$ is an iterated bundle whose base and fibers $Y$ are all of
the following sort. There are specific vertices $u,w$ with $\f(u)=w$. The set
$Z_w \subset \Chat_w$ is fixed, as is the degree $p$ of $\f_u:\Chat_u \to \Chat_w$.
Furthermore, for each $\zeta \in Z_w$ the {\em ramification data} of $\f_u$
at $\zeta$ is specified; i.e.\/ the monodromy of $\f_u^{-1}$ in a small loop 
around each $\zeta$, thought of as a conjugacy class in the symmetric group on
$p$ letters. Then $Y$ is the space of normalized degree $p$ polynomials 
with the specified ramification data. We call $Y$ a {\em Hurwitz Variety},
and observe that each $X$ is an iterated bundle with total (complex) dimension
$q-1$ whose base and fibers are all Hurwitz varieties.

The generic case is that the monodromy of $\f_u^{-1}$ in a small loop around
each $\zeta \in Z_w$ is trivial; i.e.\/ that each $\zeta$ is a regular value.
In that case $Y$ is a Zariski open subset of $\C^{p-1}$. In fact, we can say
something more precise. Let $\Delta_p \subset \C^{p-1}$ be the discriminant
variety, i.e.\/ the set of normalized degree $p$ polynomials with a multiple
root. For each $\zeta \in Z_w$ let $\Delta_{p,\zeta}:=\Delta_p+\zeta$ be the {\em translate}
of $\Delta_p$ which parameterizes the set of normalized degree $p$
polynomials $f$ for which $\zeta$ is a critical value. Then 
$$Y = \C^{p-1} - \bigcup_{\zeta \in Z_w} \Delta_{p,\zeta}$$
It turns out that the topology of $Y$ depends only on the cardinality of
$Z_w$; see \cite{Calegari_sausages} Prop.~9.14. 
This is not obvious, since the $\Delta_{p,\zeta}$ are 
singular, and they do not intersect in general position.

\subsection{$K(\pi,1)$s}

For a finite set $Z \subset \C$ and degree $p$ let $Y_p(Z)$ denote the
Hurwitz variety of normalized degree $p$ polynomials for which no element of 
$Z$ is a critical value.

As we remarked already in \S~\ref{subsection:discriminant}, when 
$|Z|=1$ the space $Y_p(Z)$ is a $K(B_p,1)$ where $B_p$ denotes the braid
group on $p$ strands. Furthermore, when $p=2$ the space $Y_2(Z)$ may be
identified with $\C-Z$ in the obvious way, so that $Y_2(Z)$ is a $K(F_n,1)$
where $F_n$ is the free group on $n$ elements, and $n=|Z|$.

It turns out (\cite{Calegari_sausages} Thm.~9.17) 
that $Y_3(Z)$ is a $K(\pi,1)$ for any finite set $Z$. This is proved by
exhibiting an explicit $\text{CAT}(0)$ 2-complex with the homotopy type of each $Y_3(Z)$.
One component of $\XX_4$ is a $K(B_4,1)$ and all the others are nontrivial
iterated fibrations where the fibers are $Y_2(Z)$ or $Y_3(Z)$s. It follows
that every component of $\XX_4$ is a $K(\pi,1)$, and in fact so is
the shift locus $\SS_4$ itself (the same is true for simpler reasons of
$\SS_3$ and $\SS_2$).

One knows few example of algebraic varieties which are $K(\pi,1)$s, and fewer
methods to construct or certify them
(one of the few general methods, which applies to certain complements
of hyperplane arrangements, is due to Deligne \cite{Deligne}). 
Is $Y_p(Z)$ a $K(\pi,1)$ for all $p$ and all $Z$?

\subsection{Monodromy}\label{subsection:monodromy}

For each $p$ and $|Z|$ there is a natural representation (well-defined up to
conjugacy) $\pi_1(Y_p(Z)) \to B_{p|Z|}$ defined by the braiding of 
the $p|Z|$ points $f^{-1}(Z)$ in $\C$ as $f$ varies in $Y_p(Z)$. This map
is evidently injective when $p=2$ or when $|Z|=1$. Is it injective in any
other case? I do not know the answer even when $p=3$ and $|Z|=2$.

Here is one reason to be interested. There is a monodromy representation of
$\pi_1(\SS_q)$ into the `Cantor braid group' --- i.e.\/ the mapping class group
of a disk minus a Cantor set --- defined by the braiding of the (Cantor)
Julia set $J_f$ in $\C$ as $f$ varies in $\SS_q$. A priori this representation
lands in the mapping class group of the plane minus a Cantor set, but it
lifts canonically to the Cantor braid group (which is a central extension) because
every $f \in \SS_q$ acts in a standard way at infinity. If one forgets the
braiding and only considers the permutation action on the Cantor set itself,
the image in $\text{Aut}(\text{Cantor Set})$
is known to be precisely equal to the automorphism group of the 
full (one-sided) shift on a $q$ element alphabet, by a celebrated theorem of
Blanchard--Devaney--Keen \cite{Blanchard_Devaney_Keen}. However, this
action of $\pi_1(\SS_q)$ 
on the Cantor set alone is very far from faithful.

The automorphism group of the Cantor set is to the Cantor braid group as a finite
symmetric group is to a (finite) braid group. It is natural to ask: is the 
monodromy representation from $\pi_1(\SS_q)$ to the Cantor braid group
injective? It turns out that the restriction of the monodromy representation
to the image of $\pi_1(Y_p(Z))$ in $\pi_1(\SS_q)$ factors through the
representation to $B_{p|Z|}$. So a precondition for the 
monodromy representation to the Cantor braid group to be injective is that
each $\pi_1(Y_p(Z)) \to B_{p|Z|}$ should be injective.

When $q=2$ we have $\pi_1(\SS_2)=\Z$ and the monodromy representation is evidently
injective, since the Cantor braid group is torsion-free. 
With Yan Mary He and Juliette Bavard we have shown
that the monodromy representation is injective in degree 3 (work in progress).

\subsection{Big Mapping Class Groups}

The Cantor braid group, and the (closely related) mapping class group of 
the plane minus a Cantor set, are quintessential examples of what are
colloquially known as {\em big mapping class groups}. The study of these groups 
is an extremely active area of current research; for an excellent recent survey see
Aramayona--Vlamis \cite{Aramayona_Vlamis}. There are connections to the 
theory of finite type
mapping class groups (particularly to stability and uniformity
phenomena in such groups); to taut foliations of 3-manifolds; 
to pruning theory and the de-Carvalho--Hall 
theory of endomorphisms of planar trees;
to Artinizations of Thompson-like groups and universal algebra; etc. (see
\cite{Aramayona_Vlamis} for references).

One major goal of this theory --- largely unrealized as yet ---
is to develop new tools for applications to dynamics in 2 real and
1 complex dimension. Cantor sets appear in surfaces 
as attractors of hyperbolic systems
(e.g.\/ in Katok--Pesin theory \cite{Katok}), and
big mapping class groups (and some closely related objects) are relevant to the
study of their moduli. The paper \cite{Calegari_sausages} and the theory of
sausages is an explicit attempt to work out some of these 
connections in a particular case.

\subsection{Rays}

Let $\Gamma$ denote the mapping class group of the plane (which we identify with
$\C$) minus a Cantor set $K$. The Cantor braid group $\hat{\Gamma}$ is the 
universal central extension of $\Gamma$. Some of the tools discussed in
this paper may be used to study $\hat{\Gamma}$ and its subgroups in some
generality; for instance, components of $\EL$ are classifying spaces
for subgroups of $\hat{\Gamma}$.

The group $\Gamma$ acts in a natural way on the set $\RR$ of {\em isotopy
classes of proper simple rays} in $\C - K$ from $\infty$ to a point in $K$. Associated
to this action are two natural geometric actions of $\Gamma$:

\begin{enumerate}
\item{there is a natural circular order on $\RR$, so that $\Gamma$ acts
faithfully by order-preserving homeomorphisms on a certain completion of
$\RR$, the {\em simple circle}; see \cite{Bavard_Walker, Calegari_planar, 
Calegari_Chen}; and}
\item{the elements of $\RR$ are the vertices of a (connected) graph (the
{\em ray graph}) whose edges correspond to pairs of rays that may be realized 
disjointly; this graph is connected, has infinite diameter, and is Gromov--hyperbolic;
see \cite{Bavard, Calegari_blog}.}
\end{enumerate}

(Landing) rays are also a critical tool in complex dynamics, and in the 
picture developed in the previous two sections.
For $K$ a Cantor Julia set, nonsingular gradient flowlines of the Green's function 
extend continuously to $K$; the set of distinct isotopy classes of 
nonsingular flowlines associated to single $K$ form a 
clique in the ray graph. Because the ray graph is 
Gromov--hyperbolic, there is (up to bounded ambiguity) a canonical path
in the ray graph between any two such cliques; 
one can ask whether such paths are coarsely
realized by paths in $\SS_q$, and if so what geometric properties such
paths have, and how this geometry manifests itself in algebraic
properties of $\pi_1(\SS_q)$. For example: does $\pi_1(\SS_q)$ admit a
(bi-)automatic structure? (to make sense of this one should 
work with a locally finite groupoid presentation for $\pi_1(\SS_q)$). 
One piece of evidence
in favor of this is that $\SS_3$ (and, for trivial reasons, $\SS_2$)
is homotopy equivalent to a locally $\text{CAT}(0)$ complex, 
and it is plausible that the same holds
for all $\SS_q$. Although there are known examples of groups which
are locally $\text{CAT}(0)$ but not bi-automatic \cite{Leary_Minasyan}, 
nevertheless in practice these two properties often go hand in hand.

\subsection{Left orderability}

A group is left-orderable if it admits a total order that is preserved under
left multiplication. The left-orderability
of braid groups (see \cite{Dehornoy_Dynnikov_Rolfsen_Wiest}) 
is key to some of their most important properties (e.g.\/ faithfulness of
the Lawrence--Kraamer--Bigelow representations \cite{Bigelow}). Left-orderability
of 3-manifold groups is also conjecturally (\cite{Boyer_Gordon_Watson}) related
to both symplectic topology (via Heegaard Floer homology) and to 
big mapping class groups via the theory of taut
foliations and universal circles; see e.g.\/ 
\cite{Calegari_Dunfield, Calegari_foliations}.
The Cantor braid group is left-orderable (via the faithful action of $\Gamma$ on 
the simple circle) so to show that $\pi_1(\SS_q)$
is left-orderable it would suffice to prove injectivity of the monodromy
representation as in \S~\ref{subsection:monodromy}.

\subsection{Comparison with finite braids}

Define $Y_q:=\C^{q-1}-\Delta_q$, the space of normalized degree $q$ polynomials
without multiple roots. Our study of $\SS_q$ has been guided by a
heuristic that one should think of $\SS_q$ as a sort of `dynamical cousin'
to $Y_q$, and that they ought to share many key algebraic and geometric properties.
Table~\ref{finitetable} compares some of what is known about
the topology of $Y_q$ and $\SS_q$. 

\begin{table}[ht]
{\small
\centering
\begin{tabular}{c||c|c|c|c} 
 & $Y_q$ & $\SS_2,\SS_3$ & $\SS_4$ & $\SS_q,q>4$ \\
\hline
locally $\text{CAT}(0)$ & yes for $q\le 6$ & yes & unknown & unknown \\
$K(\pi,1)$ & yes & yes & yes & unknown \\
$H_*$ vanishes below middle dimension & yes & yes & yes & yes \\
$\pi_1$ is mapping class group & yes & yes & unknown & unknown \\
$\pi_1$ is left-orderable & yes & yes & unknown & unknown \\ 
$\pi_1$ is biautomatic & yes & yes & yes & unknown \\
\end{tabular}
}
\vskip 12pt
\caption{Comparison of $\SS_q$ with discriminant complements $Y_q$}\label{finitetable}
\end{table}


\section{Acknowledgements}
I would like to thank Lvzhou Chen, Toby Hall, Sarah Koch, Curt McMullen, Sandra Tilmon,
Alden Walker and Amie Wilkinson for useful feedback on early drafts of this paper.

\end{document}